# Computability and analysis:
# the legacy of Alan Turing

Jeremy Avigad

Departments of Philosophy and Mathematical Sciences
Carnegie Mellon University, Pittsburgh, USA
avigad@cmu.edu

Vasco Brattka

Faculty of Computer Science
Universität der Bundeswehr München, Germany and
Department of Mathematics and Applied Mathematics
University of Cape Town, South Africa and
Isaac Newton Institute for Mathematical Sciences
Cambridge, United Kingdom
Vasco.Brattka@cca-net.de

October 31, 2012

## 1   Introduction

For most of its history, mathematics was algorithmic in nature. The geometric claims in Euclid's *Elements* fall into two distinct categories: "problems," which assert that a construction can be carried out to meet a given specification, and "theorems," which assert that some property holds of a particular geometric configuration. For example, Proposition 10 of Book I reads "To bisect a given straight line." Euclid's "proof" gives the construction, and ends with the (Greek equivalent of) Q.E.F., for *quod erat faciendum*, or "that which was to be done." Proofs of theorems, in contrast, end with Q.E.D., for *quod erat demonstrandum*, or "that which was to be shown"; but even these typically involve the construction of auxiliary geometric objects in order to verify the claim.

Similarly, algebra was devoted to developing algorithms for solving equations. This outlook characterized the subject from its origins in ancient Egypt and Babylon, through the ninth century work of al-Khwarizmi, to the solutions to the quadratic and cubic equations in Cardano's *Ars Magna* of 1545, and to Lagrange's study of the quintic in his *Réflexions sur la résolution algébrique des équations* of 1770.

The theory of probability, which was born in an exchange of letters between Blaise Pascal and Pierre de Fermat in 1654 and developed further by Christian

Huygens and Jakob Bernoulli, provided methods for calculating odds related to games of chance. Abraham de Moivre's 1718 monograph on the subject was entitled *The Doctrine of Chances: or, a Method for Calculating the Probabilities of Events in Play*. Pierre de Laplace's monumental *Théorie analytique des probabilités* expanded the scope of the subject dramatically, addressing statistical problems related to everything from astronomical measurement to the measurement in the social sciences and the reliability of testimony. Even so, the emphasis remained fixed on explicit calculation.

Analysis had an algorithmic flavor as well. In the early seventeenth century, Cavalieri, Fermat, Pascal, and Wallis developed methods of computing "quadratures," or areas of regions bounded by curves, as well as volumes. In the hands of Newton, the calculus became a method of explaining and predicting the motion of heavenly and sublunary objects. Euler's *Introductio in Analysis Infinitorum* of 1748 was the first work to base the calculus explicitly on the notion of a *function*; but, for Euler, functions were given by piecewise analytic expressions, and once again, his focus was on methods of calculation.

All this is not to say that all the functions and operations considered by mathematicians were computable in the modern sense. Some of Euclid's constructions involve a case split on whether two points are equal or not, and, similarly, Euler's piecewise analytic functions were not always continuous. In contrast, we will see below that functions on the reals that are computable in the modern sense are necessarily continuous. And even though Euler's work is sensitive to the rates of convergence of analytic expressions, these rates were not made explicit. But these are quibbles, and, generally speaking, mathematical arguments through the eighteenth century provided informal algorithms for finding objects asserted to exist.

The situation changed dramatically in the nineteenth century. Galois' theory of equations implicitly assumed that all the roots of a polynomial exist *somewhere*, but Gauss' 1799 proof of the fundamental theorem of algebra, for example, did not show how to compute them. In 1837, Dirichlet considered the example of a "function" from the real numbers to the real numbers which is equal to 1 on the rationals and 0 on the irrationals, without pausing to consider whether such a function is calculable in any sense. The Bolzano-Weierstraß Theorem, first proved by Bolzano in 1817, asserts that any bounded sequence of real numbers has a convergent subsequence; in general, there will be no way of computing such a subsequence. Riemann's proof of the open mapping theorem is based on the Dirichlet principle, an existence principle that is not computationally valid. Cantor's work on the convergence of Fourier series led him to consider transfinite iterations of point-set operations, and, ultimately, to develop the abstract notion of set.

Although the tensions between conceptual and computational points of view were most salient in analysis, other branches of mathematics were not immune. For example, in 1871, Richard Dedekind defined the modern notion of an *ideal* in a ring of algebraic integers, and defined operations on ideals in a purely extensional way. In other words, the operations were defined in such a way that they do not presuppose any particular representation of the ideals, and definitions do not indicate how to compute the operations in terms of such representations. More dramatically, in 1890, Hilbert proved what is now known as the Hilbert Basis Theorem. This asserts that, given any sequence $f_1, f_2, f_3, \ldots$ of multivariate polynomials over a Noetherian ring, there is some $n$ such that for every $m \geq n$, $f_m$ is in the ideal generated by $f_1, \ldots, f_n$. Such an $m$ cannot be computed by surveying elements of the sequence, since it is not even a continuous function on the space of sequences; even if a sequence $x, x, x, \ldots$ starts out looking like a constant sequence, one cannot rule out the possibility that the element 1 will eventually appear.

Such shifts were controversial, and raised questions as to whether the new, abstract, set-theoretic methods were appropriate to mathematics. Set-theoretic paradoxes in the early twentieth century raised the additional question as to whether they are even consistent. Brouwer's attempt, in the 1910's, to found mathematics on an "intuitionistic" conception raised a further challenge to modern methods, and in 1921, Hermann Weyl, Hilbert's best student, announced that he was joining the Brouwerian revolution. The twentieth century *Grundlagenstreit*, or "crisis of foundations," was born.

At that point, two radically different paths were open to the mathematical community:

- Restrict the methods of mathematics so that mathematical theorems have direct computational validity. In particular, restrict methods so that sets and functions asserted to exist are computable, as well as infinitary mathematical objects and structures more generally; and also ensure that quantifier dependences are also constructive, so that a "forall-exists" statement asserts the existence of a computable transformation.

- Expand the methods of mathematics to allow idealized and abstract operations on infinite objects and structures, without concern as to how these objects are represented, and without concern as to whether the operations have a direct computational interpretation.

Mainstream contemporary mathematics has chosen decisively in favor of the latter. But computation is important to mathematics, and faced with a nonconstructive development, there are options available to those specifically interested in computation. For example, one can look for computationally valid versions of nonconstructive mathematical theorems, as one does in computable and computational mathematics, constructive mathematics, and numerical analysis. There are, in addition, various ways of measuring the extent to which ordinary theorems fail to be computable, and characterizing the data needed to make them so.

With Turing's analysis of computability, we now have precise ways of saying what it means for various types of mathematical objects to be computable, stating mathematical theorems in computational terms, and specifying the data relative to which operations of interest are computable. Section 2 thus discusses *computable analysis*, whereby mathematical theorems are made computationally significant by stating the computational content explicitly.

There are still communities of mathematicians, however, who are committed to developing mathematics in such a way that every concept and assertion has an *implicit* computational meaning. Turing's analysis of computability is useful here, too, in a different way: by representing such styles of mathematics in formal axiomatic terms, we can make this implicit computational interpretation mathematically explicit. Section 3 thus discusses different styles of constructive mathematics, and the computational semantics thereof.

## 2   Computable analysis

### 2.1   From Leibniz to Turing

An interest in the nature of computation can be found in the seventeenth century work of Leibniz (see [47]). For example, his *stepped reckoner* improved on earlier mechanical calculating devices like the one of Pascal. It was the first calculating machine that was able to perform all four basic arithmetical operations, and it earned Leibniz an external membership of the British Royal Society at the age of 24. Leibniz's development of calculus is better known, as is the corresponding priority dispute with Newton. Leibniz paid considerable attention to choosing notations and symbols carefully in order to facilitate calculation, and his use of the integral symbol $\int$ and the d symbol for derivatives have survived to the present day. Leibniz's work on the binary number system, long before the advent of digital computers, is also worth mentioning. A more important contribution to the study of computation was his notion of a *calculus ratiocinator*, that is, a calculus of reasoning. Such a calculus, Leibniz held, would allow one to resolve disputes in a purely mathematical fashion:[1]

> The only way to rectify our reasonings is to make them as tangible as those of the Mathematicians, so that we can find our error at a glance, and when there are disputes among persons, we can simply say: Let us calculate, without further ado, to see who is right.

His attempts to develop such a calculus amount to an early form of symbolic logic.

With this perspective, it is not farfetched to see Leibniz as initiating a series of developments that culminate in Turing's work. Norbert Wiener has described the relationship in the following way [232]:

> The history of the modern computing machine goes back to Leibniz and Pascal. Indeed, the general idea of a computing machine is nothing but a mechanization of Leibniz's *calculus ratiocinator*. It is, therefore, not at all remarkable that the theory of the present computing machine has come to meet the later developments of the algebra of logic anticipated by Leibniz. Turing has even suggested that the problem of decision, for any mathematical situation, can always be reduced to the construction of an appropriate computing machine.

---

[1]Leibniz, *The Art of Discovery*, 1685 [233].

## 2.2   From Borel to Turing

Perhaps the first serious attempt to express the mathematical concepts of a computable real number and a computable function on the real numbers were made by Émil Borel around 1912, the year that Alan Turing was born. Borel defined computable real numbers as follows:[2]

> We say that a number $\alpha$ is computable if, given a natural number $n$, we can obtain a rational number that differs from $\alpha$ by at most $\frac{1}{n}$.

Of course, before the advent of Turing machines or any other formal notion of computability, the meaning of the phrase "we can obtain" remained vague. But Borel provided the following additional information in a footnote to that phrase:

> I intentionally leave aside the practical length of operations, which can be shorter or longer; the essential point is that each operation can be executed in finite time with a safe method that is unambiguous.

This makes it clear that Borel had an intuitive notion of an algorithm in mind. Borel then indicated the importance of number representations, and argued that decimal expansions have no special theoretical value, whereas continued fraction expansions are not invariant under arithmetic operations and hence of no practical value. He went on to discuss the problem of determining whether two real numbers are equal:

> The first problem in the theory of computable numbers is the problem of equality of two such numbers. If two computable numbers are unequal, this can obviously be noticed by computing both with sufficient precision, but in general it will not be known *a priori*. One can make clear progress in determining a lower bound on the difference of two computable numbers, whose definitions satisfy known conditions.

In modern terms, what Borel seems to recognize here is that although there is no algorithm that decides whether two computable real numbers are equal, the inequality relation between computable reals is, at least, computably enumerable. He then discussed a notion of the *height* of a number, which is based on counting the number of steps needed to construct that number, in a certain way. This concept can be seen as an early forerunner of the concept of Kolmogorov complexity. Borel considered ways that this concept might be utilized in addressing the equality problem.

---

[2]All citations of Borel are from [19], which is a reprint of [18]. The translations here are by the authors of this article; obvious mistakes in the original have been corrected.

In another section of his paper, Borel discussed the concept of a computable real number function, which he defined as follows:

> We say that a function is computable if its value is computable for any computable value of the variable. In other words, if $\alpha$ is a computable number, one has to know how to compute the value of $f(\alpha)$ with precision $\frac{1}{n}$ for any $n$. One should not forget that, by definition, to be given a computable number $\alpha$ just means to be given a method to obtain an arbitrary approximation to $\alpha$.

It is worth noting that Borel only demanded computability at computable inputs in his definition. His definition is vague in the sense that he did not indicate whether he had in mind an algorithm that transfers a *method* to compute $\alpha$ into a *method* to compute $f(\alpha)$ (which would later become known as *Markov computability* or the *Russian approach*) or whether he had in mind an algorithm that transfers an *approximation* of $\alpha$ into an *approximation* of $f(\alpha)$ (which is closer to what we now call a computable function on the real numbers, under the *Polish approach*). He also did not say explicitly that his algorithm to compute $f$ is meant to be uniform, but this seems to be implied by his subsequent observation:

> A function cannot be computable, if it is not continuous at all computable values of the variable.

A footnote to this observation then indicates that he had the Polish approach in mind:[3]

> In order to make the computation of a function effectively possible with a given precision, one additionally needs to know the modulus of continuity of the function, which is the [...] relation [...] between the variation of the function values with respect to the variation of the variable.

Borel went on to discuss different types of discontinuous functions, including those that we now call Borel measurable. In fact, the entire discussion of computable real numbers and computable real functions is preliminary to Borel's development of measure theory, and the discussion was meant to motivate aspects of that development.

---

[3]Hence Borel's result that computability implies continuity *cannot* be seen as an early version of the famous theorem of Ceĭtin, in contrast to what is suggested in [215]. The aforementioned theorem states that any Markov computable function is already effectively continuous on computable inputs and hence computable in Borel's sense, see sections 2.5 and 3, Figure 5. While this is a deep result, the observation that computable functions in Borel's sense are continuous is obvious.

### 2.3 Turing on computable analysis

Turing's landmark 1936 paper [208] is titled "On computable numbers, with an application to the Entscheidungsproblem." It begins as follows:

> The "computable" numbers may be described briefly as the real numbers whose expressions as a decimal are calculable by finite means. Although the subject of this paper is ostensibly the computable *numbers*, it is almost equally easy to define and investigate computable functions of an integrable variable or a real or computable variable, computable predicates, and so forth. The fundamental problems involved are, however, the same in each case, and I have chosen the computable numbers for explicit treatment as involving the least cumbrous technique. I hope shortly to give an account of the relations of the computable numbers, functions, and so forth to one another. This will include a development of the theory of functions of a real variable expressed in terms of computable numbers. According to my definition, a number is computable if its decimal can be written down by a machine.

At least two things are striking about these opening words. The first is that Turing chose not to motivate his notion of computation in terms of the ability to characterize the notion of a computable function from $\mathbb{N}$ to $\mathbb{N}$, or the notion of a computable *set* of natural numbers, as in most contemporary presentations; but, rather, in terms of the ability to characterize the notion of a computable real number. In fact, Section 10 is titled "Examples of large classes of numbers which are computable." There, he introduced the notion of a computable function on computable real numbers, and the notion of computable convergence for a sequence of computable numbers, and so on; and argued that, for example, $e$ and $\pi$ and the real zeros of the Bessel functions are computable. The second striking fact is that he also flagged his intention of developing a full-blown theory of computable real analysis. As it turns out, this was a task that ultimately fell to his successors, as we explain below.

The precise definition of a computable real number given by Turing in the original paper can be expressed as follows: a real number $r$ is *computable* if there is a computable sequence of 0's and 1's with the property that the fractional part of $r$ is equal to the real number obtained by prefixing that sequence with a binary point. There is a slight problem with this definition, however, which Turing discussed in a later correction [209], published in 1937. Suppose we have a procedure that, for every $i$, outputs a rational number $q_i$ with the property that

$$|r - q_i| < 2^{-i}.$$

Intuitively, in that case, we would also want to consider $r$ to be a computable real number, because we can compute it to any desired accuracy. In fact, it is not hard to show that this second definition coincides with the first: a real number has a computable binary expansion if and only if it is possible to compute it in the second sense. In other words, the two definitions are extensionally equivalent.

The problem, however, is that it is not possible to pass *uniformly* between these two representations, in a computable way. For example, suppose a procedure of the second type begins to output the sequence of approximations $1/2, 1/2, 1/2, \ldots$. Then it is difficult to determine what the first binary digit is, because at some point the output could jump just above or just below $1/2$. This intuition can be made precise: allowing the sequence in general to depend on the halting behavior of a Turing machine, one can show that there is no algorithmic procedure which, given a description of a Turing machine describing a real number by a sequence of rational approximations, computes the digits of $r$. On the other hand, for any fixed description of the first sort, there is a computable description of the second sort: either $r$ is a dyadic rational, which is to say, it has a finite binary expansion; or waiting long enough will always provide enough information to determine the digits. So the difference only shows up when one wishes to talk about computations which take descriptions of real numbers as input. In that case, as Turing noted, the second type of definition is more natural; and, as we will see below, these are the *descriptions* of the computable reals that form the basis for computable analysis.

Turing's second representation of computable reals, presented in his correction [209], is given by the formula

$$(2i-1)n + \sum_{r=1}^{\infty} (2c_r - 1) \left(\frac{2}{3}\right)^r,$$

where $i$ and $n$ provide the integer part of the represented number and the binary sequence $c_r$ the fractional part. This representation is essentially what has later been called a *signed-digit representation* with base $\frac{2}{3}$. It is interesting that Turing acknowledged Brouwer's influence (see also [72]):

> This use of overlapping intervals for the definition of real numbers is due originally to Brouwer.

In the 1936 paper, in addition to discussing individual computable real numbers, Turing also defined the notion of a computable function on real numbers. Like Borel, he adopted the standpoint that the input to such a function needs to be computable itself:

> We cannot define general computable functions of a real variable, since there is no general method of describing a real number, but we can define a computable function of a computable variable.

A few years later Turing introduced oracle machines [211], which would have allowed him to handle computable functions on *arbitrary* real inputs, simply by considering the input as an oracle given from outside and not as being itself computed in some specific way. But the 1936 definition ran as follows. First, Turing extended his definition of computable real numbers $\gamma_n$ from the unit interval to all real numbers using the formula $\alpha_n = \tan(\pi(\gamma_n - \frac{1}{2}))$. He went on:

> Now let $\varphi(n)$ be a computable function which can be shown to be such that for any satisfactory[4] argument its value is satisfactory. Then the function $f$, defined by $f(\alpha_n) = \alpha_{\varphi(n)}$, is a computable function and all computable functions of a computable variable are expressible in this form.

Hence, Turing's definition of a computable function on the real numbers essentially coincides with Markov's definition, which would later form the basis of the Russian approach to computable analysis (see sections 2.6 and 3.5). Computability of $f$ is not defined in terms of approximations, but in terms of functions $\varphi$ that transfer algorithms that describe the input into algorithms that describe the output. There are at least three objections against the technical details of Turing's definition, some of which have already been addressed by Guido Gherardi in [72], which is the first careful discussion of Turing's contributions to computable analysis:

1. Turing used the binary representation of real numbers in order to define $\gamma_n$, and hence the derived notion of a computable real function is somewhat unnatural and restrictive. However, we can consider this problem as being repaired by Turing's second approach to computable reals in [209], which we have described above.

2. Turing should have mentioned that the function $f$ is only well-defined by $f(\alpha_n) = \alpha_{\varphi(n)}$ if $\varphi$ is extensional in the sense that $\alpha_n = \alpha_k$ implies $\alpha_{\varphi(n)} = \alpha_{\varphi(k)}$. However, we can safely assume that this is what he had in mind, since otherwise his equation does not even define a single-valued function.

3. Turing did not allow arbitrary suitable computable functions $\varphi$, but he restricted his definition to total functions $\varphi$. However, it is an easy exercise to show that for any partial computable $\varphi$ that maps all satisfactory numbers to satisfactory numbers there is a total computable function $\psi$ such that $\psi(n)$ and $\varphi(n)$ are numbers of the same total computable function for all satisfactory inputs $n$. Hence the restriction to total $\varphi$ is not an actual restriction.

These three considerations support our claim that Turing's definition of a computable real function is essentially the same as Markov's definition.

It is worth pointing out that Turing introduced another important theme in computable analysis, namely, developing computable versions of theorems of analysis. For instance, he pointed out:

---

[4] Turing calls a natural number $n$ *satisfactory* if, in modern terms, $n$ is a Gödel index of a total computable function.

... we cannot say that a computable bounded increasing sequence of computable numbers has a computable limit.

This shows that Turing was aware of the fact that the Monotone Convergence Theorem cannot be proved computably. He did, however, provide a (weak) computable version of the Intermediate Value Theorem:

(vi) If $\alpha$ and $\beta$ are computable and $\alpha < \beta$ and $\varphi(\alpha) < 0 < \varphi(\beta)$, where $\varphi(a)$ is a computable increasing continuous function, then there is a unique computable number $\gamma$, satisfying $\alpha < \gamma < \beta$ and $\varphi(\gamma) = 0$.

The fact that the limit of a computable sequence of computable real numbers need not be computable motivates the introduction of the concept of computable convergence:

We shall say that a sequence $\beta_n$ of computable numbers *converges computably* if there is a computable integral valued function $N(\varepsilon)$ of the computable variable $\varepsilon$, such that we can show that, if $\varepsilon > 0$ and $n > N(\varepsilon)$ and $m > N(\varepsilon)$, then $|\beta_n - \beta_m| < \varepsilon$. We can then show that

(vii) A power series whose coefficients form a computable sequence of computable numbers is computably convergent at all computable points in the interior of its interval of convergence.

(viii) The limit of a computably convergent sequence is computable.

And with the obvious definition of "uniformly computably convergent":

(ix) The limit of a uniformly computably convergent computable sequence of computable functions is a computable function. Hence

(x) The sum of a power series whose coefficients form a computable sequence is a computable function in the interior of its interval of convergence.

From (viii) and $\pi = 4(1 - \frac{1}{3} + \frac{1}{5} - ...)$ we deduce that $\pi$ is computable. From $e = 1 + 1 + \frac{1}{2!} + \frac{1}{3!} + ...$ we deduce that $e$ is computable. From (vi) we deduce that all real algebraic numbers are computable. From (vi) and (x) we deduce that the real zeros of the Bessel functions are computable.

Here we need to warn the reader that Turing erred in (x): the mere computability of the sequences of coefficients of a power series does not guarantee computable convergence on the whole interior of its interval of convergence, but only on each compact subinterval thereof (see Theorem 7 in [174] and Exercise 6.5.2 in [224]).

Before we leave Turing behind, some of his other, related work is worth mentioning. For example, in 1938, Turing published an article, "Finite Approximations to Lie Groups" [210], in the *Annals of Mathematics*, addressing the question as to which infinite Lie groups can be approximated by finite groups. The idea of approximating infinitary mathematical objects by finite ones is central to computable analysis. The problem that Turing considered here is more restrictive in that he required that the approximants be groups themselves, and the paper is not in any way directly related to computability; but the parallel is interesting. What is more directly relevant is the fact that Turing was not only interested in computing real numbers and functions in theory, but in practice. The Turing archive[5] contains a sketch of a proposal, in 1939, to build an analog computer that would calculate approximate values for the Riemann zeta-function on the critical line. His ingenious method was published in 1943 [212], and in a paper published in 1953 [214] he described calculations that were carried out in 1950 on the Manchester University Mark 1 Electronic Computer. Finally, we mention that, in a paper of 1948 [213], Turing studied linear equations from the perspective of algebraic complexity. There he introduced the concepts of a *condition number* and of *ill-conditioned* equations, which are widely used in numerical mathematics nowadays. (For a discussion of this work, see Lenore Blum's contribution to this collection, as well as [72].) In some very interesting unpublished work Turing provided an algorithm to compute an absolutely normal real number (see the discussion by Becher et. al. in [13]) and in this way he even made an early contribution to algorithmic randomness.

## 2.4   From Turing to Specker

Turing's ideas on computable real numbers where taken up in 1949 by Ernst Specker [201], who completed his Ph.D. in Zürich under Hopf in the same year. Specker was probably attracted to logic by Paul Bernays, whom he thanks in that paper. Specker constructed a computable monotone bounded sequence $(x_n)$ whose limit is not a computable real number, thus establishing Turing's claim. In modern terms one can easily describe such a sequence: given an injective computable enumeration $f : \mathbb{N} \to \mathbb{N}$ of a computably enumerable but non-computable set $K \subseteq \mathbb{N}$, the partial sums

$$x_n = \sum_{i=0}^{n} 2^{-f(i)}$$

form a computable sequence without a computable limit. Such sequences are now sometimes called *Specker sequences*.

It is less well-known that Specker also studied different representations of real numbers in that paper, including the representation by computably converging Cauchy sequences, the decimal representation, and the Dedekind cut representation (via characteristic functions). In particular, he considered both primitive recursive and general computable representations in all three instances.

---

[5] See AMT/C/2 in `http://www.turingarchive.org/browse.php/C`.

His main results include the fact that the class of real numbers with primitive recursive Dedekind cuts is strictly included in the class of real numbers with primitive recursive decimal expansions, and the fact that the latter class is strictly included in the class of real numbers with primitive recursive Cauchy approximation. He also studied arithmetic properties of these classes of real numbers.

Soon after Specker published his paper, Rosza Péter included Specker's discussion of the different classes of primitive recursive real numbers in her book, *Rekursive Funktionen* [173]. In a review of that book [184], Raphael M. Robinson noted that all the representations of real numbers mentioned above yield the same class of numbers if one uses computability in place of primitive recursiveness. However, this is no longer true if one moves from single real numbers to sequences of real numbers. Mostowski [161] proved that for computable sequences of real numbers the analogous strict inclusions hold, as they were proved by Specker for single primitive recursive real numbers: the class of sequences of real numbers with uniformly decidable Dedekind cuts is strictly included in the class of sequences of real numbers with uniformly computable decimal expansion, which in turn is strictly included in the class of sequences of real numbers with computable Cauchy approximations. Mostowski pointed out that the striking coincidence between the behavior of classes of single primitive recursive real numbers and computable sequences of real numbers is not very well understood:

> It remains an open problem whether this is a coincidence or a special case of a general phenomenon whose cause ought to be discovered.

A few years after Specker and Robinson, the logician H. Gordon Rice [183] proved what Borel had already observed informally, namely that equality is not decidable for computable real numbers. In contrast, whether $a < b$ or $a > b$ holds can be decided for two non-equal computable real numbers $a, b$. Rice also proved that the set $\mathbb{R}_c$ of computable real numbers forms a real algebraically closed field, and that the Bolzano-Weierstraß Theorem does not hold computably, in the following sense: there is a bounded c.e. set of computable numbers without a computable accumulation point. Kreisel wrote a review of Rice's paper for the *Mathematical Reviews* and pointed out that he had already proved results that are more general than Rice's; namely, he proved that any power series with a computable sequence of coefficients (that are not all identical to zero) has only computable zeros in its interval of convergence [121], and that there exists a computable bounded set of rational numbers which contains no computable subsequence with a computable modulus of convergence [120].

In a second paper on computable analysis [202], Specker constructed a computable real function that attains its maximum in the unit interval at a non-computable number. In a footnote he mentioned that this solves an open problem posed by Grzegorczyk [82], a problem that Lacombe had already solved independently [137]. Specker mentioned that a similar result can be derived from a theorem of Zaslavskiĭ [239]. He also mentioned that the Intermediate Value Theorem has a (non-uniform) computable version, but that it fails for sequences of functions.

### 2.5 Computable analysis in Poland and France

In the 1950's, computable analysis received a substantial number of contributions from Andrzej Grzegorczyk [82, 83, 84, 85] in Poland and Daniel Lacombe [82, 83, 84, 85, 124, 122, 123, 136, 137, 138, 139, 140, 141, 142, 143] in France. Grzegorczyk is the best known representative of the Polish school of computable analysis, a school of research that came into being soon after Turing published his seminal paper. Apparently, members of the Polish school of functional analysis became interested in questions of computability, and Stefan Banach and Stanisław Mazur gave a seminar talk on this topic in Lviv (now in Ukraine, at that time Polish) on January 23, 1937 [11]. Unfortunately, the second world war got in the way, and Mazur's course notes on computable analysis (edited by Grzegorczyk and Rasiowa) where not published until much later [154].

Banach and Mazur's notion of computability for functions on the real numbers was defined using computable sequences. Accordingly, we call a function $f : \mathbb{R}_c \to \mathbb{R}_c$ *Banach-Mazur computable* if it maps computable sequences to computable sequences. While this concept works well for well-behaved functions (see, for example, the discussion of Pour-El and Richards' work on linear functions in Section 2.6), it yields a weaker concept of computability of real number functions in general, as proved later by Aberth (again, see Section 2.6). One of Mazur's theorems was that every Banach-Mazur computable function is already continuous on the computable reals. This theorem can be seen as a forerunner and even as a stronger version of a theorem of Kreisel, Lacombe, and Shoenfield [123], which was independently proved in slightly different versions by Ceïtin [36] and later by Moschovakis [160]. Essentially, the theorem says that every real number function that is computable in the sense of Markov is continuous (on the computable reals). (In this particular case one even has the stronger conclusion that the function is effectively continuous. See Figure 5 in Section 2.9.)

Grzegorczyk's work [84] was mostly concerned with computable functions $f : [0,1] \to \mathbb{R}$ on all real numbers of the unit interval and he proved a number of equivalent characterizations of this concept. For instance $f$ is computable if and only if each of the following conditions hold:

1. $f$ maps computable sequences to computable sequences and it has a computable modulus of uniform convergence,

2. $f$ can be represented by a computable function $h$ on rational intervals that is monotone with respect to inclusion and that approximates $f$ in the sense that $\{f(x)\} = \bigcap\{h(I) : x \in I\}$ for all $x \in [0,1]$.

Later on, the first characterization became the basic definition of computability used by Pour-El and Richards (see section 2.6), whereas the second characterization is the starting point for domain theoretic characterizations of computability [196, 54].

Lacombe's pioneering work opened up two research directions in computable analysis which were to play an important role. First, he initiated the study of

c.e. open sets of reals, which are open sets that can be effectively obtained as unions of rational intervals [139]. Together with Kreisel, Lacombe proved that there are c.e. open sets that contain all computable reals and that have arbitrarily small Lebesgue measure [124]. Secondly, Lacombe anticipated the study of computability on more abstract spaces than Euclidean space, and was one of the first to define the notion of a computable metric space [143].

## 2.6  Computable analysis in North America

More comprehensive approaches to computable analysis were developed in the United States during the 1970's and 1980's. One stream of results is due to Oliver Aberth, and is accumulated in his book, *Computable Analysis* [2]. His work is built on Markov's notion of computability for real number functions. We call a function $f : \mathbb{R}_c \to \mathbb{R}_c$ on computable reals *Markov computable* if there is an algorithm that transforms algorithms that describe inputs $x$ into algorithms that describe outputs $f(x)$. Among other results, Aberth provided an example of a Markov computable function $f : \mathbb{R}_c \to \mathbb{R}_c$ that cannot be extended to a computable function on all real numbers [2, Theorem 7.3]. Such a function $f$ can, for instance, be constructed with a Specker sequence $(x_n)$, which is a computable monotone increasing bounded sequence of real numbers that converges to a non-computable real $x$. Using effectively inseparable sets, such a Specker sequence can even be constructed so that it is effectively bounded away from any computable point (see, for instance, the Effective Modulus Lemma in [176]). One can use such a special Specker sequence and place increasingly steeper peaks on each $x_n$, as in Figure 1. If adjusted appropriately, this construction leads to a Markov computable function $f : \mathbb{R}_c \to \mathbb{R}_c$ that cannot be extended to a computable function on the whole of $\mathbb{R}$. In fact, it cannot even be extended to a continuous function on $\mathbb{R}$, since the peaks accumulate at $x$. Pour-El and Richards [176, Theorem 6, page 67] later refined this counterexample to obtain a differentiable, non-computable function $f : \mathbb{R} \to \mathbb{R}$ whose restriction to the computable reals $\mathbb{R}_c$ is Markov computable. This can be achieved by replacing the peaks in Aberth's construction by smoother peaks of decreasing height.

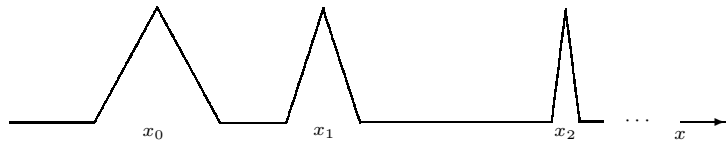

Figure 1: Construction of a function $f : \mathbb{R} \to \mathbb{R}$

Another well-known counterexample due to Aberth [1] is the construction of a function $F : \mathbb{R}_c^2 \to \mathbb{R}_c$ that is Markov computable and uniformly continuous on a rectangle centered around the origin, with the property that the differential equation

$$y'(x) = f(x, y(x))$$

with initial condition $y(0) = 0$ has no computable solution $y$, not even on an arbitrarily small interval around the origin. This result can be extended to computable functions $f : \mathbb{R}^2 \to \mathbb{R}$, which was proved somewhat later by Pour-El and Richards [177]. It can be interpreted as saying that the Peano Existence Theorem does not hold computably.

Marian Pour-El and J. Ian Richards' approach to computable analysis was the dominant approach in the 1980's, and is nicely presented in their book [176]. The basic idea is to axiomatize computability structures on Banach spaces using computable sequences and a characterization of computability that is derived from Grzegorczyk's first characterization, which was mentioned in section 2.5. This approach is tailor-made for well-behaved maps such as linear ones. One of their central results is their First Main Theorem [179], which states that a linear map $T : X \to Y$ on computable Banach spaces is computable if and only if it is continuous and has a c.e. closed graph, i.e.

$$T \text{ computable } \iff T \text{ continuous and } \mathrm{graph}(T) \text{ c.e. closed.}$$

Moreover, the theorem states that a linear function $T$ with a c.e. closed graph satisfies

$$T \text{ computable } \iff T \text{ maps computable inputs to computable outputs.}$$

Here a closed set is called *c.e. closed* if it contains a computable dense sequence. The crucial assumption that the $\mathrm{graph}(T)$ is c.e. closed has sometimes been ignored tacitly in presentations of this theorem and this has created the wrong impression that the theorem identifies computability with continuity. In any case, the theorem is an important tool and provides many interesting counterexamples. This is because the contrapositive version of the theorem implies that every linear discontinuous $T$ with a c.e. closed graph maps some computable input to a non-computable output. For instance, the operator of differentiation

$$d :\subseteq C[0,1] \to C[0,1], f \mapsto f'$$

is known to be linear and discontinuous and it is easily seen to have a c.e. closed graph (for instance the rational polynomials are easily differentiated) and hence it follows that there is a computable and continuously differentiable real number function $f : [0,1] \to \mathbb{R}$ whose derivative $f'$ is not computable. This is a fact that has also been proved directly by Myhill [163].

The First Main Theorem is applicable to many other operators that arise naturally in physics, such as the wave operator [178, 175]. In particular, it follows that there is a computable wave $(x, y, z) \mapsto u(x, y, z, 0)$ that transform into a non-computable wave $(x, y, z) \mapsto u(x, y, z, 1)$ if it evolves for one time unit according to the wave equation

$$\frac{\partial^2 u}{\partial x^2} + \frac{\partial^2 u}{\partial y^2} + \frac{\partial^2 u}{\partial z^2} - \frac{\partial^2 u}{\partial t^2} = 0$$

from initial condition $\frac{\partial u}{\partial t} = 0$. This result has occasionally been interpreted as suggesting that it might be possible to build a wave computer that would violate Church's Thesis. However, a physically more appropriate analysis shows that this is not plausible (see Section 2.7 and also [172]). Pour-El and Richards' definitions and concepts regarding Banach spaces have been transferred to the more general setting of metric spaces by Yasugi, Mori and Tsujii [235] in Japan.

Another stream of research on computable analysis in the United States came from Anil Nerode and his students and collaborators. For instance, Metakides, Nerode and Shore proved that the Hahn-Banach Theorem is not computably valid, but admits a non-uniform computable solution in the finite dimensional case [157, 158]. Joseph S. Miller studied sets that appear as the set of fixed points of a computable self map of the unit ball in Euclidean space [159]. Among other things he proved that a closed set $A \subseteq [0,1]^n$ is the set of fixed points of a computable function $f : [0,1]^n \to [0,1]^n$ if and only if it is co-c.e. closed and has a co-c.e. closed connectedness component. Zhou also studied computable open and closed sets on Euclidean space [242]. In a series of papers [103, 104, 102], Kalantari and Welch provided results that shed further light on the relation between Markov computable functions and computable functions. Douglas Cenzer and Jeffrey B. Remmel have studied the complexity of index sets for problems in analysis [37, 38, 39]. McNicholl and his co-authors have started a study of computable complex analysis [156, 155, 152, 153].

It is possible to study computable functions of analysis from the point of view of computational complexity. We cannot discuss such results in detail here, but in passing we mention the important work of Harvey Friedman and Ker-I Ko [113, 114, 115, 117, 64, 118], which is well summarized in [116]; and, more recently, interesting work of Mark Braverman, Stephen Cook and Akitoshi Kawamura [28, 27, 29, 105, 106]. Braverman and Michael Yampolsky have recently written a nice monograph on computability properties of Julia sets [30].

## 2.7 Computable analysis in Germany

In Eastern Germany, Jürgen Hauck produced a sequence of interesting papers on computable analysis that were only published in German [86, 87, 88], building on the work of Dieter Klaua [107, 108]. The novel aspect of Hauck's work is that he studied the notion of a representation (of real numbers or other objects) in its own right. This aspect has been taken up in the work of Weihrauch's school of computable analysis [224], which has its origins in the work of Christoph Kreitz and Klaus Weihrauch [125, 126, 226, 222]. Kreitz and Weihrauch began to develop a systematic theory of representations that is centered around the notion of an *admissible* representation. Representations $\delta_X :\subseteq \mathbb{N}^{\mathbb{N}} \to X$ are (potentially partial) maps that are surjective. Admissible representations $\delta_X$ are particularly well-behaved with respect to some topology given on the represented set $X$. If $\delta_X$ and $\delta_Y$ are admissible representations of topological spaces $X$ and $Y$, respectively, then a (partial) function $f :\subseteq X \to Y$ is continuous if and only if there exists a continuous $F :\subseteq \mathbb{N}^{\mathbb{N}} \to \mathbb{N}^{\mathbb{N}}$ on Baire space $\mathbb{N}^{\mathbb{N}}$ such that the diagram in Figure 2 commutes.

The perspective arising from the concept of a representation as a map sheds new light on representations of real numbers as they were studied earlier by Specker, Mostowski and others. Their results can be developed in a more uniform way, using the concept of computable reducibility for representations, which yields a lattice of real number representations. The position of a representation in this lattice characterizes the finitely accessible information content of the respective representation, see Figure 3. In particular, the behavior of representations that was mysterious to Mostowski can be explained naturally in this way.

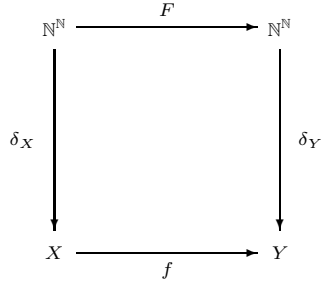

Figure 2: Notions of computable real number functions

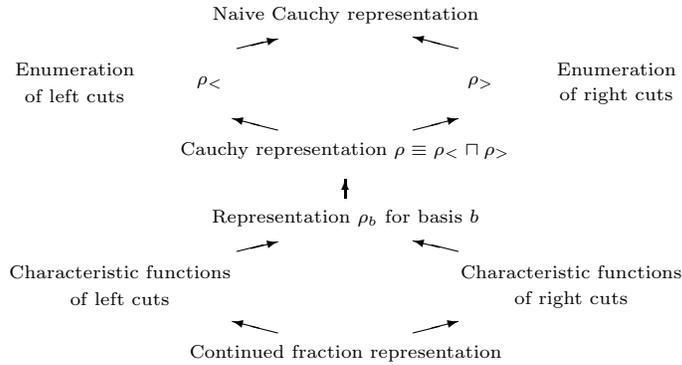

Figure 3: The lattice of real number representations

The concept of an admissible representation has been extended to a larger classes of topological spaces by the work of Matthias Schröder [193]. In his definition, an admissible representation with respect to a certain topology on the represented space is just one that is maximal among all continuous representations with respect to computable reducibility. Hence, admissibility can be seen as a completeness property.

Many other interesting results have been developed using the concept of a representation (see, for instance, the tutorial [25] for a more detailed survey). For instance, Peter Hertling [91] proved that the real number structure is computably categorical in the sense that up to computable equivalence there is one and only one representation that makes all the ingredients of the structure computable. He also studied a computable version of Riemann's Mapping Theorem [90], and proved that the Hyperbolicity Conjecture implies that the Mandelbrot set $M \subseteq \mathbb{R}^2$ is computable in the sense that its distance function is computable. Finally, he also proved that there is a Banach-Mazur computable function $f : \mathbb{R}_c \to \mathbb{R}_c$ that is not Markov computable [92].

Matthias Schröder and Klaus Weihrauch have studied a very general setting for computational complexity in analysis using proper representations [190, 225, 194] and Schröder has also studied online complexity of real number functions [191, 192], following earlier work of Schönhage [189]. Robert Rettinger and Klaus Weihrauch studied the computational complexity of Julia sets [180]. In another interesting stream of papers weaker notions of computability on real numbers were studied [240, 181, 182] by Xizhong Zheng, Klaus Weihrauch and Robert Rettinger.

Ning Zhong and Klaus Weihrauch studied computability of distributions [241] and proved computability of the solution operator of several partial differential equations [227, 228, 229, 230, 231]. In particular, they revised the wave equation and proved that if this equation is treated in physically realistic spaces, then it turns out to be computable [227].

Martin Ziegler has studied models of hypercomputation and related notions of computability for discontinuous functions [243]. Martin Ziegler and Stéphane Le Roux have revisited co-c.e. closed sets. Among other things, they proved that the Bolzano-Weierstraß Theorem is not limit computable, in the sense that there is a computable sequence in the unit cube without a limit computable cluster point [144]. Baigger [10] proved that the Brouwer Fixed Point Theorem has a computable counterexample (essentially following Orevkov's ideas [167], who proved the corresponding result for Markov computability much earlier).

Just as the computational content of representations of real numbers can be analyzed in a lattice, there is another lattice that allows one to classify the computational content of (forall-exists) theorems. This lattice is based on a computable reduction for functions that has been introduced by Klaus Weihrauch and that has recently been studied more intensively by Arno Pauly, Guido Gherardi, Alberto Marcone, Brattka, and others [171, 23, 21]. Theorems such as the Hahn-Banach Theorem [73], the Baire Category Theorem, the Intermediate Value Theorem [22], the Bolzano-Weierstraß Theorem [24] and the Radon-Nikodym Theorem [100] have been classified in this lattice as well as finding Nash equilibria or solving linear equations [170]. Figure 4 shows the respective relative positions of these theorems in the Weihrauch lattice. This classification yields a uniform perspective on the computational content of all these theorems, and the position of a theorem in this lattice fully characterizes certain properties of the respective theorem regarding computability. In particular, it implies all aforementioned computability properties of these theorems.

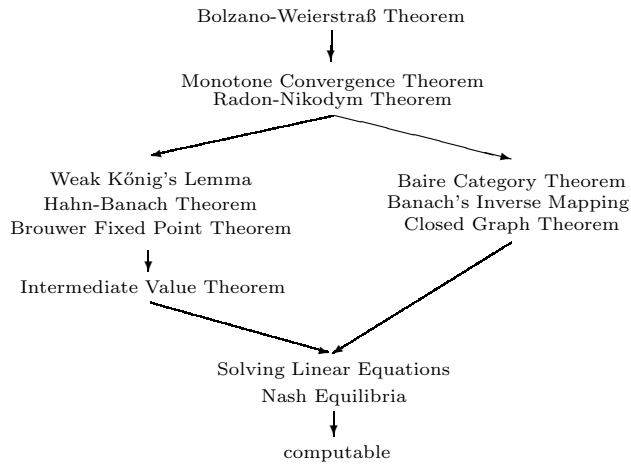

Figure 4: The computational content of theorems in the Weihrauch lattice

## 2.8 Computable measure theory

The field of *algorithmic randomness* relies on effective notions of measure and measurable set. Historically, debates in the 1930's as to how to make sense of the notion of a "random" sequence of coin flips gave rise to the Kolmogorov axioms for a probability space, and the measure-theoretic interpretation of probability. On that interpretation, the statement that a "random" element of a space has a certain property should be interpreted as the statement that the set of counterexamples has measure 0. In that framework, talk of random elements is only a manner of speaking; there is no truly "random" sequence of coin flips $x$, since each such element is contained in the null set $\{x\}$.

In the 1960's, Per Martin-Löf, a student of Kolmogorov, revived the notion of a random sequence of coin flips [149, 150]. The idea is that, from a computability-theoretic perspective, it makes sense to call a sequence random if it avoids all suitably effective null sets. For example, a *Martin-Löf test* is an effectively null $G_\delta$ set, which is to say, an effective sequence of computably open sets $G_i$ with the property that for each $i$, the Lebesgue measure of $G_i$ is less than $2^{-i}$. An sequence of coin flips *fails* such a test if it is in the intersection of the sets $G_i$. An element of the unit interval is said to be *Martin-Löf random* if it does not fail any such test. Osvald Demuth [48] defined an equivalent notion, independently, in 1975.

Martin-Löf randomness can be given alternate characterizations in terms of computable betting strategies, or in terms of information content. There are by now many other notions of randomness in the literature. One, based on the influential work of C.-P. Schnorr [188, 187], is similar to Martin-Löf randomness but requires that, moreover, the measures of the open sets $G_i$ are uniformly computable. It is impossible to even begin to survey this very active area here, but detailed overviews can be found in the monographs [165, 50]. Kučera and Nies [135] provides a nice discussion of Demuth's contributions to the subject.

Measure theory, more generally, has been studied from a computability-theoretic perspective. Most work in algorithmic randomness has focused on random elements of $2^{\mathbb{N}}$, or, equivalently, random elements of the unit interval $[0, 1]$. But there are a number of ways of developing the notion of a computable measure on more general spaces; see the work of Gács [66], Hertling and Weihrauch [93] and Hoyrup and Rojas [98]. Weihrauch [223] and Schröder [195] have discussed different representations of probability measures and Weihrauch and Wu [234] have further developed the foundations of computable measure theory. Edalat [53] has studied computable measure and integration theory from the perspective of domain theory. Bosserhoff [20] and Hoyrup and Rojas [98] have defined various notions of computability for measurable functions. This makes it possible to consider results in measure theory, probability, and dynamical systems in computability-theoretic terms.

For example, computable aspects of the ergodic theorems have been studied extensively. The pointwise and mean ergodic theorems assert that if $T$ is a measure-preserving transformation of a finite measure space and $f$ is any measurable function, then the ergodic averages $A_n f = \frac{1}{n} \sum_{i<n} f \circ T^i$ converge pointwise, and in the $L^2$ norm, respectively. V'yugin [220, 221] and, independently, Avigad and Simic [9, 197, 5] have shown that the mean and pointwise ergodic theorems do not hold computably, but Avigad, Gerhardy, and Townsner [8] have shown that a rate of convergence can be computed from $T$, $f$, and the norm of the limit. In particular, the rate of convergence can be computed from $T$ and $f$ when the transformation is ergodic. V'yugin [220, 221] has shown that if $T$ is a measurable transformation of the unit interval and $f$ is a computable function, then the sequence of ergodic averages $(A_n f)$ converges on every Martin-Löf random real. Franklin and Towsner [60] have recently shown that this is sharp. Hoyrup and Rojas [69] have generalized V'yugin's result to computable measurable functions. Indeed, the work of Gács, Galatolo, Hoyrup, and Rojas [68, 67], together with the results of [8], shows that the Schnorr random reals are exactly the ones at which the ergodic theorem holds with respect to every ergodic measure and computable measure-preserving transformation (see also the discussion in [166]). Bienvenu, Day, Hoyrup, and Shen [15] and, independently, Franklin, Greenberg, Miller, and Ng [59] have extended V'yugin's result to the characteristic function of a computably open set, in the case where the measure is ergodic.

Other theorems of measure-theory and measure-theoretic probability have been considered in this vein. (Such work is, to some extent, foreshadowed by the study of measure theory in the context of reverse mathematics, by Simpson, Yu, and others; see, for example, [238, 236, 237].) The Lebesgue differentiation theorem states, roughly speaking, that for almost every point of $\mathbb{R}^n$, the value of an integrable function at that point is equal to the limit of increasingly small averages taken around that point. Pathak [169] showed that, when the integrable function is computable in an appropriate sense, the theorem is true at every Martin-Löf random point. Pathak, Rojas, and Simpson [166] later strengthened the conclusion to the Schnorr-random real numbers, and showed that this is sharp; this same result was obtained independently by Rute [186]. Another theorem of Lebesgue's states that a function of bounded variation is differentiable almost everywhere. Recasting results of Demuth [49], Brattka, Miller, and Nies [26] have shown that when the function in question is computable, the result holds at Martin-Löf random real numbers, and that this is sharp. Moreover, they obtained clean characterizations of other notions of randomness in terms of differentiability. Freer, Kjos-Hansen, Nies, and Stephan [61] have obtained similar results for the class of Lipschitz functions. Rute [186] has considered a number of differentiability and martingale convergence theorems, characterizing the points at which the theorems hold, and determining when rates of convergence are and are not computable. Ackerman, Freer, and Roy have studied the computability of conditional distributions [3], Hoyrup, Rojas, and Weihrauch have studied the computability of Radon-Nikodym derivatives [99], and Freer and Roy have provided a computable version of de Finetti's theorem [62].

## 2.9 Appendix: Notions of computability for real number functions

In this appendix, we summarize the various notions of computability for real number functions that were discussed in previous sections, and we briefly indicate the logical relations between these notions. For simplicity we restrict the discussion to functions that are total on all real numbers or on all computable real numbers.

A *rapidly converging Cauchy name* for a real number $x$ is a sequence $(x_i)_{i \in \mathbb{N}}$ of rationals that converges to $x$ rapidly, i.e. in such a way that for every $i$ and $j \geq i$, $|x_j - x_i| < 2^{-i}$. (Any other computable rate of convergence would work equally well.)

1. $f : \mathbb{R} \to \mathbb{R}$ is called *computable* if there is an algorithm that transforms each given rapidly converging Cauchy name of an arbitrary $x$ into such a name for $f(x)$ (on a Turing machine with one-way output tape).

2. $f : \mathbb{R}_c \to \mathbb{R}_c$ is called *Borel computable* if $f$ is computable in the previously described sense, but restricted to computable reals.

3. $f : \mathbb{R}_c \to \mathbb{R}_c$ is called *Markov computable* if there is an algorithm that converts each given algorithm for a computable real $x$ into an algorithm for $f(x)$. Here an algorithm for a computable real number $x$ produces a rapidly converging Cauchy name for $x$ as output.

4. $f : \mathbb{R}_c \to \mathbb{R}_c$ is called *Banach-Mazur computable* if $f$ maps any given computable sequence $(x_n)$ of real numbers into a computable sequence $(f(x_n))$ of real numbers.

The diagram in Figure 5 summarizes the logical relations between different notions of computability. Some of these equivalences can be generalized to partial functions with well-behaved domains, and to other computable metric spaces (see [89, 92]).

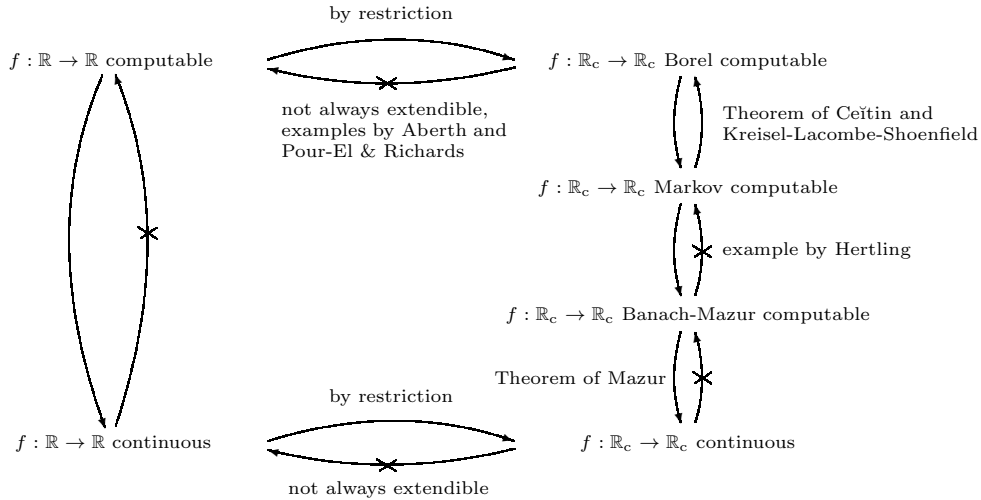

Figure 5: Notions of computable real number functions

## 3   Constructive analysis

As noted in Section 1, the growing divergence of conceptual and computational concerns in the early twentieth century left the mathematical community with two choices: either restrict mathematical language and methods to preserve direct computational meaning, or allow more expansive language and methods and find a place for computational mathematics within it. We have seen how the theory of computability has supported the latter option: given broad notions of real number, function, space, and so on, we can now say, in precise mathematical terms, which real numbers, functions, and spaces are *computable*.

There is still, however, a community of mathematicians that favors the first option, namely, adopting a "constructive" style of mathematics, whereby theorems maintain a direct computational interpretation. This simplifies language: rather than say "there exists a computable function $f$" one can say "there exists a function $f$," with the general understanding that, from a constructive standpoint, proving the existence of such a function requires providing an algorithm to compute it. Similarly, rather than say "there is an algorithm which, given $x$, computes $y$ such that . . . ," one can say "given $x$, there exists $y$ such that . . . " with the general understanding that a constructive proof, by its very nature, provides such an algorithm.

In sum, from a constructive point of view, all existence claims are expected to have computational significance. We will see, however, that this general stricture is subject to interpretation, especially when it comes to infinitary mathematical objects. Today there are many different "flavors" of constructive mathematics, with varying norms and computational interpretations. Here the methods of logic and computability theory prove to be useful yet again: once we design a formal deductive system that reasonably captures a style of constructive reasoning, we can provide a precise computational semantics. This serves to clarify the sense in which a style of constructive mathematics has computational significance, and helps bridge the two approaches to preserving computational meaning.

This section will provide a brief historical overview of a some of the different approaches to constructive mathematics, as well as some associated formal axiomatic systems and their computational semantics. Troelstra and van Dalen [207] and Beeson [14] provide more thorough introductions to these topics.

## 3.1 Kroneckerian mathematics

In 1887, facing the rising Dedekind-Cantor revolution, the mathematician Leopold Kronecker wrote an essay in which he urged the mathematical community to maintain focus on symbolic and algorithmic aspects of mathematics.

> [A]ll the results of the profoundest mathematical research must in the end be expressible in the simple forms of the properties of integers. But to let these forms appear simply, one needs above all a suitable, surveyable manner of expression and representation for the numbers themselves. The human spirit has been working on this project persistently and laboriously since the greyest prehistory.... ([130, page 955])

One gets a better sense of Kronecker's views by considering the way they played out in his work, such as his treatment of algebraic numbers in his landmark *Grundzüge einer arithmetischen Theorie der algebraischen Grössen* [128], or in his treatment of the fundamental theorem of algebra in *Ein Fundamentalsatz der allgemeinen Arithmetik* [129]. In general, Kronecker avoided speaking of "arbitrary" real numbers and functions, and, rather, dealt with algebraic *systems* of such things, wherein the objects are given by symbolic expressions, and operations on the objects are described in terms of symbolic calculations. Of course, it is still an important task to understand how rational approximations can be obtained from an algebraic description of a real number. But even when dealing with convergent sequences and limits (as in, for example, the proof of Dirichlet's theorem on primes in an arithmetic progression in his lectures on number theory), Kronecker always provided explicit information as to how quickly a sequence converges. The approach is nicely described by his student Hensel:

> He believed that one can and must in this domain formulate each definition in such a way that its applicability to a given quantity can be assessed by means of a finite number of tests. Likewise that an existence proof for a quantity is to be regarded as entirely rigorous only if it contains a method by which that quantity can really be found. ([131], quoted in Stein [203, p. 250].)

A number of articles by Harold Edwards (including [55, 56, 57]) provide an excellent overview of Kronecker's outlook and contributions to mathematics.

## 3.2 Brouwerian intuitionism

During the 1910's, the Dutch mathematician L. E. J. Brouwer advanced a new philosophy of mathematics, *intuitionism*, which made strong pronouncements as to the appropriate methods of mathematical reasoning. Brouwer held that mathematical objects are constructions of the mind, and that we come to know a mathematical theorem only by carrying out a mental construction that enables us to see that it is true. In particular, seeing that a statement of the form $A \wedge B$ is true involves seeing that $A$ is true, and that $B$ is true; seeing that $A \to B$ is true involves having a mental procedure that transforms any construction witnessing the truth of $A$ to a construction witnessing the truth of $B$; seeing that $\neg A$ is true involves having a procedure that transforms a construction witnessing the truth of $A$ to a contradiction (that is, something which, evidently, cannot be the case); and seeing that $A \vee B$ is true requires carrying out a mental construction enabling one to see that $A$ is true, or carrying out a mental construction to see that $B$ is true; and similarly for statements involving the universal and existential quantifiers. This account, developed further by Heyting and Kolmogorov, has come to be known as the BHK interpretation. It is not hard to see that, according to this interpretation, the law of the excluded middle, $A \vee \neg A$, does not generally hold. For example, if $A$ is the statement of the Goldbach conjecture, then we cannot presently assert $A \vee \neg A$, because we do not currently know that the Goldbach conjecture is true, nor do we know that it is false.

There is a strong solopsistic strand to Brouwer's philosophy, in that the mathematical knowledge one has is a reflection of one's inner mental life, independent of an external world or other thinkers. He also held that mathematical knowledge and thought are independent of language, which he took to be a deeply limited and flawed means of communicating mathematical ideas. Intuitionism may seem to have little to do with computation *per se*, but, as we will see below, if one replaces Brouwer's mental constructions and procedures with symbolic representations and algorithms, one is left with an essentially computational view of mathematics. Thus, Brouwer's mathematical ideas have proved to be influential in computer science, despite the different motivations.

Brouwer was also an influential topologist, and an important part of his intuitionist program was to develop an intuitionistically appropriate foundation for reasoning about the continuum [32, 33]. To that end, he introduced the notion of a *choice sequence*, that is, an "arbitrary" sequence of natural numbers. Some choice sequences are generated by a law, such as the sequence $0, 0, 0, \ldots$. At the other extreme, some sequences are "lawless," that is, generated by events that are not predetermined. Brouwer introduced the *continuity principle*, which, roughly speaking, asserts that any function from choice sequences to the natural numbers is continuous; in other words, the value of the function at a given choice sequence depends only on a finite initial segment of that sequence. (See Troelstra and van Dalen [207, I.4.6] for a precise formulation.)

Brouwer went on to develop an intuitionistic set theory based on such sequences. A *spread*, in his terminology, is a tree on $\mathbb{N}$ such that each node has at least one successor. Brouwer saw this data as a way of specifying a collection of choice sequences; that is, a choice sequence is in the spread if and only if every initial segment is in the tree. A *fan* is a spread with the property that every node has only finitely many successors. At the risk of clouding some intuitions, we can translate these notions to classical terms: a choice sequence is an arbitrary element of Baire space (the space of functions from $\mathbb{N}$ to $\mathbb{N}$ under the product topology), a spread corresponds to a closed subset of Baire space, and a fan corresponds to a compact subset of Baire space.

In 1927, Brouwer [34] introduced the *bar theorem*, which, despite the name, is essentially an axiomatic principle that provides transfinite induction on well-founded trees on $\mathbb{N}$. Roughly speaking, it asserts that any property that holds outside and at the leaves of such a well-founded tree, and is preserved in passing from all the children of a node to the node itself, holds of every finite sequence. An immediate corollary is the *fan theorem*, which asserts that a well-founded fan is finite. Using the fan theorem and the continuity principle, Brouwer showed that every function from $[0, 1]$ to $\mathbb{R}$ is uniformly continuous. When properties of choice sequences are understood in terms of computable predicates on Baire space, this formulation of analysis accords quite well with Grzegorczyk's notion of computability, as discussed in Section 2.5.

For more information on intuitionistic mathematics, see [52, 94, 207, 31]. For more on Brouwer and his philosophical views, see [216, 217].

## 3.3 Early formal systems for constructive mathematics

Now let us consider some formal axiomatic systems that capture such constructive styles of reasoning. A theory now known as *primitive recursive arithmetic* was designed by Thoralf Skolem [199] to capture Hilbert's notion of "finitary" mathematical reasoning. In its strictest form, the theory has variables ranging over the natural numbers, but no quantifiers. One starts with some basic functions on the natural numbers, and is allowed to define new functions using a schema of *primitive recursion*:

- $f(0, z_1, \ldots, z_n) = g(z_1, \ldots, z_n)$, and

- $f(x + 1, z_1, \ldots, z_n) = h(x, f(x, z_1, \ldots, z_n), z_1, \ldots, z_n)$,

where $g$ and $h$ have previously been defined. Although this schema may seem limited, with ingenuity one can show that almost any reasonable function (and hence relation) on the natural numbers is primitive recursive. In particular, the primitive recursive relations are closed under bounded quantification. As a result, primitive recursive arithmetic can be viewed as a reasonable framework in which to carry out a Kroneckarian constructivism.

Indeed, primitive recursive arithmetic is often taken as a starting point for formalizing basic mathematical notions in many presentations of formal axiomatic foundations for mathematics, classical or constructive. It plays a central role in Hilbert and Bernays' landmark *Grundlagen der Mathematik* [96], in Kleene's *Introduction to Metamathematics* [111], and Goodstein's *Recursive Number Theory* [79]. Goodstein's later *Recursive Analysis* [80] develops topics like integration and differentiation on the basis of primitive recursive arithmetic.

Today, finitism is viewed as an extreme form of intuitionism. Foundational writings in the 1920's, however, do not show a clear distinction between the two. Indeed, as late as 1930, von Neumann used the two terms interchangeably in his exposition of formalism [219] at the Second Conference for Epistemology and the Exact Sciences in Königsberg. (This was, incidentally, the meeting where Gödel announced the first incompleteness theorem.)

The situation was clarified when Arend Heyting, a student of Brouwer's, provided formal axiomatizations of intuitionistic mathematics, despite Brouwer's antipathy towards formalization. Specifically, he provided formal treatments of intuitionistic logic, arithmetic, and Brouwer's theory of choice sequences, which appeared in a series of three papers published in 1930 [95]. In 1933, Gödel [76], and later Gentzen [70], showed that classical first-order arithmetic could be interpreted in intuitionistic first-order arithmetic via what is now known as the *double-negation translation*. At the time, classical first-order arithmetic, which extends primitive recursive arithmetic with quantifiers ranging over the natural numbers and induction for all formulas in the language, was viewed as strictly stronger than finitism. The interpretation of classical arithmetic in intuitionistic arithmetic thereby prompted the realization that intuitionism and finitism diverge as well.

### 3.4 Realizability

The developments described in this section so far all predate Turing's analysis of computability. Once notions of computability were in place, however, it was not long before Stephen Kleene showed that they can be used to provide a precise sense in which intuitionistic mathematics has a computational interpretation. Towards the end of a paper [109] published in 1943, but first presented to the American Mathematical Society in 1940, he put forth the following:

> Thesis III. A proposition of the form $(x)(Ey)A(x, y)$ containing no free variables is provable constructively, only if there is a general recursive function $\phi(x)$ such that $(x)A(x, \phi(x))$. ([109, Section 16, page 69])

One can turn this thesis into *theorems* by showing that the statement holds when we replace "constructive provability" by provability in particular axiomatic systems. Kleene carried out this program over the next few years, together with his student, David Nelson [110, 164]. For that purpose, Kleene introduced the notion of "realizability," which is, roughly speaking, an analysis of the Brouwer-Heyting-Kolmogorov interpretation in computability-theoretic terms. Specifically, Kleene defined a relation *e realizes* $\varphi$ inductively, where $e$ is a natural number and $\varphi$ is a sentence in the language of first-order arithmetic, with the following clauses:

- If $\varphi$ is atomic, then $e$ realizes $\varphi$ if and only if it is true.

- $e$ realizes $\theta \wedge \eta$ if and only if $e$ is of the form $2^a 3^b$, where $a$ realizes $\theta$ and $b$ realizes $\eta$.

- $e$ realizes $\theta \vee \eta$ if and only if $e$ is of the form $2^0 3^a$ and $a$ realizes $\theta$, or $e$ is of the form $2^1 3^b$ and $b$ realizes $\eta$.

- $e$ realizes $\theta \rightarrow \eta$ if and only if $e$ is the Gödel number of a partial recursive function $f$ which, given any $a$ realizing $\theta$, returns a number $f(a)$ realizing $\eta$.

- $e$ realizes $\exists x \, \theta(x)$ if and only if $e$ is of the form $2^x 3^a$ and $a$ realizes $\theta(\bar{x})$, where $\bar{x}$ is the numeral that denotes $x$.

- $e$ realizes $\forall x \, \theta(x)$ if and only if $e$ is the Gödel number of a partial recursive function $f$ which, given any $x$, returns a number $f(x)$ realizing $\theta(\bar{x})$.

More generally, if $\varphi$ has free variables, $e$ realizes $\varphi$ if and only if it realizes its universal closure. A formula $\varphi$ is said to be *realizable* if there is some number realizing it.

Kleene emphasized that this does not provide a reductive analysis of intuitionistic truth, insofar as quantifiers and logical connectives themselves appear in the definition. In particular, if $\varphi$ is a true, purely universal sentence, then anything realizes $\varphi$; and a realizer for a negated sentence carries no useful information at all. However, for any sentence $\varphi$ in the language of arithmetic, the statement "$e$ realizes $\varphi$" can also be expressed in the language of arithmetic, and one can inquire as to how these two assertions are related to one another. Nelson showed that if intuitionistic logic proves $\varphi$, then it proves that $\varphi$ is realizable. Unfortunately, it is not the case that for any $\varphi$, intuitionistic arithmetic proves that $\varphi$ is equivalent to its own realizability; but Nelson showed that intuitionistic arithmetic *does* prove that this equivalence is realizable. Moreover, one can strengthen the clause for implication:

- $e$ realizes $\theta \rightarrow \eta$ if and only if $\theta$ implies $\eta$ *and* $e$ is the Gödel number of a partial recursive function $f$ which, given any $a$ realizing $\theta$, returns a number $f(a)$ realizing $\eta$.

In that case, intuitionistic arithmetic proves that a formula $\varphi$ is true if and only if it is realizable.

The idea that a constructive proof provides explicit "evidence" for the truth of a theorem in question, and that such evidence can often be described in computability-theoretic terms, is a powerful one. It can help illuminate the "meaning" or "computational content" of a formal axiomatic system; and, as a purely technical device, it can be used to obtain metamathematical properties of such systems, such as provability and unprovability results. In 1967, Kleene and Vesley [112] presented a formalization of Brouwerian analysis together with a suitable realizability interpretation.

By now, there is a bewildering array of realizability relations in the literature. Variations can depend on any of the following features:

- the language expressing the mathematical assertions (first-order, second-order, or higher-order, etc.);

- the kinds of realizers (arbitrary computable functions, computable functions in a particular class, or a broader class of functions);

- the descriptions of the realizers (e.g. whether they are represented by natural numbers, or terms in a formal language);

- whether or not the realizability relation itself is expressed in a formal system; or

- the particular clauses of the realizability relation itself (see, for example, the variant of the clause for implication above).

All these decisions have bearing on the axioms and rules that are realized, and the metamathematical consequences once can draw. See Troelstra [206] for a definitive reference, as well as [14, 205] for more information. Realizability theory can also be used to translate results from constructive analysis into computable analysis; see [12, 145].

## 3.5 The Russian school of constructive mathematics

The post-Turing era brought a new approach to constructive mathematics, the principal tenets of which were set forth by A. A. Markov in the late 1940's and developed through the 1950's [147, 148]. The result is what has come to be known as the "Markov school" or "Russian school" of constructive mathematics, with contributions by Nikolai Shanin, I. Zaslavskiĭ, Gregory Ceĭtin, Osvald Demuth, and Boris Kushner, among many others.[6] Aspects of their work have already been discussed above in connection with computable analysis. Indeed, a hallmark of this style of constructivity is that it relies explicitly on notions of computability, which is to say, the real numbers are explicitly defined to be computable reals; the notion of a function from the natural numbers to the natural numbers is explicitly defined to be a computable function; and so on. In contrast to contemporary computable analysis, however, the Russian school insisted that proofs also have a constructive character, so that, for example, a proof of a statement of the form $\forall x \; \exists y \; \varphi(x, y)$ statement can be seen to yield a computable dependence of $y$ on $x$. This style of constructivity can be viewed as "constructive computable mathematics," or "constructive recursive mathematics," and, indeed, is often referred to as such.

---

[6] A detailed bibliography can be found on the *Computability and Complexity in Analysis Network*, http://cca-net.de/publications/.

Note that this style of constructivity stands in stark contrast to Brouwerian intuitionism: even if Brouwer could have identified his "constructions" with formal notions of computability, he would have been unlikely to do so. The Russian school also adopted a principle that is not found in Brouwerian intuitionism, namely, Markov's principle. This states that if $P$ is a decidable property of natural numbers (that is, $P(x) \vee \neg P(x)$ holds for every $x$), and it is contradictory that no $x$ satisfies $P$, then some $x$ satisfies $P$. In symbols:

$$\neg \forall x \, \neg P(x) \to \exists x P(x).$$

The intuition is that one can find an $x$ satisfying $P(x)$ by simply searching for it systematically, since the hypothesis guarantees that the search cannot fail to turn up such an $x$.

One can find good overviews of this style of constructivity in books by Aberth [2] and Kushner [133], as well as Kushner's survey [134].

## 3.6 Theories of finite types

In constructive recursive mathematics, one can interpret functions as Turing-machine indices, which can, in turn, be thought of as descriptions of a computer program. Thus, in constructive recursive mathematics, sets and functions can be "coded," or represented, by natural numbers. But ordinary mathematics deals not just with sets and functions of natural numbers, but also with sets of sets, sets of functions, functionals defined on spaces of functions, and so on. Rather than represent all of these using indices for computable objects, it is more natural, for some purposes, to keep the computational interpretation implicit, and take mathematical objects at face value.

To that end, it is convenient to adopt the language of *finite types*. The idea traces back to foundational frameworks designed by Frege, Russell and Whitehead, and Church, which will be discussed in Section 3.8. Roughly, a "type" can be thought of as a syntactic classification of mathematical objects. To obtain the finite types, start with the type $\mathsf{N}$, intended to denote the natural numbers, so that an object of type $\mathsf{N}$ is a natural number. Add the rule that whenever $\mathsf{A}$ and $\mathsf{B}$ are finite types, so is $\mathsf{A} \to \mathsf{B}$. Intuitively, an object of type $\mathsf{A} \to \mathsf{B}$ is a function from $\mathsf{A}$ to $\mathsf{B}$. Thus, we can form the type $\mathsf{N} \to \mathsf{N}$ of functions from $\mathsf{N}$ to $\mathsf{N}$, the type $(\mathsf{N} \to \mathsf{N}) \to \mathsf{N}$ of functionals from $\mathsf{N} \to \mathsf{N}$ to $\mathsf{N}$, and so on. It is sometimes also convenient to add a base type $\mathsf{Bool}$ for the Boolean values "true" and "false," and product types $\mathsf{A} \times \mathsf{B}$, but these are inessential.

Following Gödel [77], we can extend the set of primitive recursive functions to the set of *primitive recursive functionals of finite type* by extending the schema of primitive recursion to the higher types. This provides a syntactic calculus for defining objects of the various types, and Gödel's theory $T$ provides a calculus for reasoning about these objects.

There are then two ways of giving $T$ a computational interpretation. The first is to remain at the level of syntax: one provides an explicit procedure to "reduce" any term in the calculus to a canonical normal form (see Tait [204] and Girard [75]). Thus if $F$ is a term of type $\mathsf{A} \to \mathsf{B}$, one can view $F$ as denoting the function which takes any term $t$, in normal form, and returns the normal form corresponding to $F(t)$.

A second approach, however, provides a more natural computational understanding of the finite types. For each type $\sigma$, define the set of *hereditarily recursive functions of type $\sigma$*, $\mathsf{HRO}_\sigma$, inductively as follows: the hereditarily recursive functions of type $\mathsf{N}$ are simply the natural numbers, and the hereditarily recursive functions of type $\mathsf{A}$ to $\mathsf{B}$ are those indices $e$ such that for every hereditarily recursive function $x$ of type $\mathsf{A}$, $\varphi_e(x)$ is defined, and is a hereditarily recursive function of type $\mathsf{B}$. Thus the hereditarily recursive functions of type $\mathsf{N} \to \mathsf{N}$ are (indices of) total computable functions; hereditarily recursive functions of type $(\mathsf{N} \to \mathsf{N}) \to \mathsf{N}$ are computable functions which, for each total computable function, return a number; and so on. It is then easy to interpret each term of $T$ as a hereditarily recursive functional.

In general, hereditarily recursive functions are not extensional: because functions act on indices, it can happen that two indices $e$ and $e'$ compute the same function from $\mathsf{N}$ to $\mathsf{N}$, and yet $F(e) \neq F(e')$ for some hereditarily recursive functional $F$. One can repair this by inductively defining extensional equality between functionals, and insisting that the interpretations of the function types preserve this equality. The resulting set of functionals is then called the *hereditarily effective operations*, $\mathsf{HEO}$. This is an instance of what computer scientists refer to as a *PER model*, since for each $\sigma$, $\mathsf{HEO}_\sigma$ is given by a partial equivalence relation, which is to say, an equivalence relation on a subset of the natural numbers.

One obtains the finite type extensions $\mathsf{HA}^\omega$ of Heyting arithmetic by extending $T$ with quantifiers and induction over arbitrary finite types. Gödel's *Dialectica interpretation* provides an interpretation of $\mathsf{HA}^\omega$ in $T$. Alternatively, one can show that provable formulas in $\mathsf{HA}^\omega$ are realized by terms in $T$. Thus, whenever $\mathsf{HA}^\omega$ proves $\forall x \, \exists y \, \varphi(x, y)$, where $x$ and $y$ are variables of any finite type, there is a hereditarily effective operation which computes $y$ from $x$. The models $\mathsf{HRO}$ and $\mathsf{HEO}$ can be formalized in intuitionistic arithmetic, $\mathsf{HA}$, showing that $\mathsf{HA}^\omega$ is also interpretable in first-order intuitionistic arithmetic. (For discussions of $\mathsf{HA}^\omega$, $\mathsf{T}$, and models thereof, see [7, 206, 207].)

### 3.7 Bishop-style constructive mathematics

In 1967, the American analyst, Errett Bishop, published a book, *Foundations of Constructive Analysis* [16], which launched a new era in constructivity. In the preface, he wrote:

> It appears ... that there are certain mathematical statements that are merely evocative, which make assertions without empirical validity. There are also mathematical statements of immediate empirical validity, which say that certain performable operations will produce certain observable results, for instance, the theorem that every pos-itive integer is the sum of four squares. Mathematics is a mixture of the real and the ideal, sometimes one, sometimes the other, often so presented that it is hard to tell which is which. The realistic component of mathematics—the desire for pragmatic interpretation—supplies the control which determines the course of development and keeps mathematics from lapsing into meaningless formalism. The idealistic component permits simplifications and opens possibilities which would otherwise be closed. The methods of proof and objects of investigation have been idealized to form a game, but the actual conduct of the game is ultimately motivated by pragmatic considerations.

For 50 years now there have been no significant changes in the rules of this game. Mathematicians unanimously agree on how mathematics should be played...

This book is a piece of constructivist propaganda, designed to show that there does exist a satisfactory alternative.

This was a landmark in the history of constructive mathematics. Brouwerian intuitionism relied not only on a vocabulary that is foreign to most working mathematicians, but also on principles, such as the continuity of every function defined on $[0, 1]$, that are not classically valid. Similarly, the Russian school's explicit restriction to computable objects sets it apart from contemporary mathematics. A central feature of Bishop's constructive mathematics is that the theorems are classically valid, and, indeed, look like ordinary mathematical theorems. The point, however, is that they are established in such a way that every theorem has computational significance. This is achieved by stating definitions and theorems carefully, restricting generality ("avoiding pseudogenerality," in Bishop's words), and adhering to methods that retain computational meaning.

*Foundations* thus began with an informal statement of constructive set-theoretic principles. To start with, "a *sequence* is a rule which associates to each positive integer $n$ a mathematical object $a_n$." Then:

The totality of all mathematical objects constructed in accord with certain requirements is called a *set.* The requirements of the construction, which vary with the set under consideration, determine that set. Thus the integers are a set, the rational numbers are a set, and the collection of all sequences each of whose terms is an integer is a set. Each set $A$ will be endowed with a relation = of equality. This relation is a matter of convention, except that it must be an *equivalence relation.*

And, as far as functions are concerned:

The dependence of one quantity on another is expressed in the basic notion of an operation. An *operation* from a set $A$ to a set $B$ is a rule $f$ which assigns an element $f(a)$ of $B$ to each element $a$ of $A$. The rule must afford an explicit, finite, mechanical reduction of the procedure for constructing $f(a)$ to the procedure for constructing
$a$... The most important case occurs when

$$f(a_1) = f(a_2)$$

whenever $a_1$ and $a_2$ are equal elements of $A$. Such an operation $f$ is
called a *function*.

Although Bishop did not provide a formal axiomatic foundation, it is reasonable to seek a formal interpretation of this framework. One option is to use a
system like $\mathsf{HA}^\omega$, and interpret each set as a definable subset of a type with a
definable equivalence relation; that is, take the membership relation $x \in A$ to
be given by a formula $\varphi_A(x)$ and take equality $x =_A y$ to be given by another
formula $\psi_{=_A}(x, y)$. The interpretations of $\mathsf{HA}^\omega$ discussed in the last section
then give Bishop-style constructive mathematics a direct computational interpretation. Bishop-style mathematics can also be developed in constructive type
theory (which will be discussed below) viewing sets as types equipped with an
equivalence relation (a.k.a. "setoids").

Bishop-style constructive mathematics is, in a sense, the most "pure" (or, at
least, minimal) constructive approach discussed so far, in that it does not rely
on bar induction or the continuity principle from intuitionistic mathematics,
nor Markov's principle from constructive recursive mathematics. This makes
the framework more appealing to classical mathematicians. For example, we
have seen that the Russian school defines a real number to be a computable real
number, and a function from reals to reals to be computable in an appropriate
sense. This means that classical mathematicians cannot view their theorems
about real numbers as such; they have to be interpreted as saying something
more restrictive. Similarly, the fact that, in a Brouwerian framework, all functions on the real numbers are continuous shows that the Brouwerian notion of
function departs from the classical one. In contrast, even though it is consistent with Bishop-style mathematics to think of all real numbers and functions
as being computable, the framework is fully consistent with ordinary classical
mathematics. In other words, the theorems in the framework can be viewed
as contemporary mathematical theorems, proved in such a way that the results
have additional computational significance.

Bishop's work was continued by a number of people, including Douglas
Bridges, Ray Mines, Fred Richman, and many others, and remains a mainstay
of modern constructivity (see, for example, [17, 31]).

## 3.8   Intuitionistic higher-order arithmetic

The finite types, discussed in Section 3.6, trace back to Gottlob's Frege foundational system [63], which was built on a single base type of individuals. They
made their way, in modified form, into the ramified type theory of Russell and
Whitehead's *Principia Mathematica* [185], and ultimately into Alonzo Church's
formulation of higher order logic as *simple type theory* [40]. All of these systems include some form of a scheme of *comprehension*,

$$\exists X \; \forall y \; (y \in X \leftrightarrow \varphi(y)),$$

which asserts that any formula $\varphi$ with a free variable $y$ of type $\sigma$ gives rise to a set, or predicate, $X$, of type $\sigma \to \mathsf{Bool}$. (More precisely, Frege's "extensions" of formulas stood as proxy for such objects.)

All of the systems just described are based on classical logic, but one can just as well consider intuitionistic versions, for example, adding comprehension axioms to $\mathsf{HA}^\omega$. The result is *intuitionistic higher-order logic*, or $\mathsf{IHOL}$. From a logical perspective, such a system is much stronger than $\mathsf{HA}^\omega$. One reason to be interested in such a system is that it represents the internal logic of a *topos*, the algebraic structure that Alexander Grothendieck used to study sheaves over a space (see, for example, [146]).

It is perhaps a matter of debate whether intuitionistic higher-order logic deserves to be called "constructive." But one thing that speaks in favor of this is that one can give $\mathsf{IHOL}$ a computational interpretation. In fact, this can be done in various ways, paralleling the various ways of giving a computational interpretation to Gödel's $T$. For example, one can define an explicit normalization procedure which reduces proofs of $\mathsf{IHOL}$ to a canonical normal form; methods based on Girard's *candidats de reducibilité* [74, 75] show the reduction procedure always terminates. Alternatively, one can give a realizability interpretation by interpreting $\mathsf{IHOL}$ in Martin Hyland's *effective topos* [101].

## 3.9 Constructive type theory

At present, the predominant foundational frameworks for constructive mathematics take the form of *constructive type theory*. Such frameworks unify two of the foundational trends we have seen so far:

- type theory, in the Frege-Russell-Church-Gödel tradition; and

- the notion of explicit "evidence" for a constructive claim, in the tradition of the BHK interpretation and realizability.

It is the use of *dependent types* that makes this unification possible.

In simple type theory, types cannot depend on parameters. For example, given a type $\mathsf{A}$ and a fixed natural number $n$, one can form the type $\mathsf{A}^n$ of $n$-tuples of elements of $\mathsf{A}$, but one cannot view these as a *family* of related types that depend on the parameter $n$. In other words, one cannot view $\mathsf{A}^n$ as a type that depends on a variable $n$ of type $\mathsf{N}$. This is exactly the sort of thing that dependent type theory is designed to support. To start with, the type $\mathsf{A} \to \mathsf{B}$ of functions which take an argument in $\mathsf{A}$ and return a value in $\mathsf{B}$ can be generalized to a dependent product $\prod_{x:\mathsf{A}} \mathsf{B}(x)$, where $\mathsf{B}(x)$ is a type that can depend on $x$. Intuitively, elements of this type are functions that map an element $a$ of $\mathsf{A}$ to an element of $\mathsf{B}(a)$. When $\mathsf{B}$ does not depend on $x$, the result is just $\mathsf{A} \to \mathsf{B}$. Similarly, product types $\mathsf{A} \times \mathsf{B}$ can be generalized to dependent sums $\sum_{x:A} B(x)$. Intuitively, elements of this type are pairs $(a, b)$, where $a$ is an element of $A$ and $b$ is an element of $B(a)$. When $B$ does not depend on $x$, this is just $A \times B$.

The second conceptual innovation in constructive type theory is to internalize the notion of constructive evidence. According to the BHK interpretation, knowing a mathematical theorem amounts to having a construction that realizes the claim. Thus we can view any mathematical proposition as specifying a special type of data, namely, the type of construction that is appropriate to realizing it. This has come to be known as the "propositions as types correspondence" or the "Curry-Howard correspondence" [45, 46, 97], developed by Curry, Howard, Tait, Martin-Löf, and Girard, among others. The point is that logical operations look a lot like operations on datatypes. For example, in propositional logic, from $A$ and $B$ one can conclude $A \wedge B$. One can read this as saying that given a proof $a$ of $A$ and a proof $b$ of $B$ of $B$ one can "pair" them to obtain a proof $(a, b)$ of $A \wedge B$. In other words, $A \wedge B$ and $A \times B$ are governed by the same rules. Similarly, the interpretation of implication $A \rightarrow B$ mirrors the rules for function types: giving a proof of $A \rightarrow B$ amounts to constructing a function from $A$ to $B$. In the same way, a proof of $\forall x : A\, B(x)$ is a function which, given any $a$ in $A$, return a proof of $B(a)$.

Thus, in constructive type theory, some types are naturally viewed as types of data and some are naturally viewed as propositions, but the two interact and are governed by similar rules. A single calculus gives the rules for defining mathematical objects and proving propositions; that is, the calculus provides a set of rules for carrying out mathematical constructions of both sorts. Two of the most commonly used systems today are *Martin-Löf type theory* [151], and the *calculus of inductive constructions* [43], which extends the original *calculus of constructions* due to Coquand and Huet [42]. The relationship between the Martin-Löf type theory and the calculus of constructions is similar to the relationship between $HA^\omega$ and intuitionistic higher-order logic: the calculus of constructions has "impredicative" comprehension principles that render it much stronger than Martin-Löf's predicative counterpart.

As was the case with $HA^\omega$ and intuitionistic higher-order logic, one can give these systems a computational interpretation in various ways. Indeed, semantics for constructive type theory draws on the full range of the theory of programming language semantics, making use of realizability interpretations, strong normalization proofs, domain theory, and more. The literature on this subject is extensive; see, for example, [4, 44].

## 3.10 Computational interpretations of classical theories

We characterized the Russian school of constructive recursive mathematics as reasoning about computable objects in a constructive way. One can maintain the first requirement axiomatically while giving up the second: for example, the subsystem of second-order arithmetic known as $RCA_0$ is a classical system for which first-order quantifiers can be interpreted as ranging over the natural numbers and second-order quantifiers can be interpreted as ranging over com-putable sets and functions. Thus $\mathsf{RCA}_0$ is a reasonable setting for formalizing classical computable analysis (see [198]).

To what extent can one preserve a computational interpretation of quantifier dependences in a formal system that includes the law of the excluded middle? There is a trivial sense in which any reasonable classical theory has a computational interpretation. Let $\varphi$ be a $\Pi_2$ statement, that is, an assertion of the form $\forall x \, \exists y \, R(x, y)$ where $x$ and $y$ range over natural numbers and $R$ is a primitive recursive relation. Suppose some classical theory $T$ proves $\varphi$. Then, assuming $T$ can be trusted in this regard, $\varphi$ is *true*, which means that a simple-minded computer program that, on input $x$, searches systematically for a $y$ satsifying $R(x, y)$ is guaranteed to succeed.

There is also a fundamental sense in which such an interpretation cannot be extended to $\Pi_3$ sentences. Let $T(e, x, s)$ be Kleene's primitive recursive relation that expresses that $s$ is a halting computation sequence for Turing machine $e$ on input $x$. Then the classically valid statement that any given Turing machine $e$ either halts on input 0 or doesn't can be expressed as follows:

$$\forall e \, \exists s \, \forall s' \, (T(e, 0, s) \lor \neg T(e, 0, s')).$$

But any function mapping an $e$ to an $s$ satisfying the conclusion provides a solution to the halting problem, and thus there is no computable function witnessing the $\forall e \, \exists s$ dependence.

Nonetheless, one can often give classical logic an *indirect* computational interpretation. One way to do this is to interpret a classical theory in a constructive one, using the double-negation translation and tricks such as the Friedman-Dragalin translation [51, 65] or the Dialectica interpretation [77, 7] to recover $\Pi_2$ theorems (see also [6, 41]). One can also provide more direct computational interpretations, such as realizability relations of various sorts. Griffin [81] has shown that classical logic can be understood in terms of a standard semantics for programming languages with control operators, such as exceptions. This computational interpretation is captured by the $\lambda\mu$-calculus designed by Parigot [168] (see also [200]). Chetan Murthy [162] has shown that this interpretation can be seen as the result of combining a double-negation translation with the Friedman-Dragalin trick, and then using intuitionistic realizability. Jean-Louis Krivine [127] has provided a realizability interpretation for full classical set theory.

Sometimes what one wants from a classical proof is not a computational interpretation *per se*, but useful computational or quantitative information that is hidden by classical methods. In the 1950's, Kreisel spoke of "unwinding" classical proofs to obtain such information, a program which has developed under the heading of "proof mining," by Kohlenbach and others [119].

*Acknowledgments.* We are grateful to Martín Escardó, Guido Gherardi, André Nies, Klaus Weihrauch, and an anonymous referee for helpful comments and suggestions. Avigad's work has been partially supported by National Science Foundation grant DMS-1068829 and Air Force Office of Scientific Research grant

FA9550-12-1-0370. Brattka's research was supported by a Marie Curie International Research Staff Exchange Scheme Fellowship within the 7th European Community Framework Programme and by the National Research Foundation of South Africa.